\documentclass[a4paper,11pt,reqno]{amsart}

\usepackage{latexsym}
\usepackage{amsthm,amssymb,amsmath,amsopn}
\usepackage{graphicx}

\newtheorem{theorem}{Theorem}[section]

\newtheorem{lemma}[theorem]{Lemma}
\newtheorem{problem}[theorem]{Problem}
\newtheorem{example}[theorem]{Example}
\newtheorem{corollary}[theorem]{Corollary}

\newtheorem{proposition}[theorem]{Proposition}

\theoremstyle{remark}

\DeclareMathOperator{\Aut}{Aut}
\DeclareMathOperator{\Fix}{Fix}
\DeclareMathOperator{\sgn}{sgn}

\newcommand{\Ind}{\big\uparrow}
\newcommand{\Res}{\big\downarrow}
\newcommand{\ind}{\!\!\uparrow}
\newcommand{\res}{\!\!\downarrow}

\DeclareMathOperator{\cAut}{\chi Aut}
\DeclareMathOperator{\clAut}{cAut}
\DeclareMathOperator{\Out}{Out}
\DeclareMathOperator{\Inn}{Inn}
\newcommand{\notdiv}{\hspace{-2.5pt}\not\hspace{2.5pt}\mid}
\begin{document}


\begin{abstract}
Let $X$ be a character table of the symmetric group  $S_n$. 
It is shown that unless $n = 4$ or $n=6$, there is a unique
way to assign partitions of~$n$ to the rows and columns of $X$ so that
for all $\lambda$ and $\nu$, $X_{\lambda\nu}$ is equal to 
$\chi^\lambda(\nu)$, 
the value of the irreducible character of $S_n$ 
labelled by $\lambda$ on elements of cycle type $\nu$. 
Analogous results are proved for alternating groups,
and for the Brauer character tables of symmetric and
alternating groups. 
\end{abstract}

\title[Labelling character tables]{Labelling the character tables of symmetric and alternating groups}


\author{Mark Wildon}
\date{October 12, 2006 } 
\email{m.j.wildon@swansea.ac.uk}

\maketitle

\section{Introduction}
In \cite{Nagao} Nagao proved that if 
$G$ is a finite group with a character table which differs from a
character table of the symmetric group~$S_n$ only by 
a permutation of its rows and columns, then~$G$ is isomorphic
to~$S_n$. In this paper we consider a question naturally
raised by Nagao's result. To state it, we must recall
that the ordinary irreducible characters of~$S_n$ are
canonically labelled by the partitions of $n$, and that
this set also labels the conjugacy classes of~$S_n$.
Given partitions $\lambda$ and $\nu$ of $n$, let $\chi^\lambda(\nu)$
be the value of the irreducible character of $S_n$ labelled
by $\lambda$ on elements of cycle type $\nu$.

Now suppose that one has 
discovered (for example, by applying Nagao's theorem) that
a given  square matrix  is an unlabelled character table of the symmetric
group~$S_n$. We ask: when can one go further, and \emph{uniquely} reconstruct
the partitions labelling its rows and columns? 
The answer is given by the following theorem.

\begin{theorem}\label{Thm:labelling}
Let $X$ be a character table of the symmetric group 
$S_n$. Unless $n=4$ or $n=6$, there is a unique
way to assign partitions of $n$ to the rows and columns of $X$ so that
$X_{\lambda\nu} = \chi^\lambda(\nu)$ for all partitions $\lambda, \nu$ of $n$.
If $n = 4$ or $n=6$ then there are exactly two different labellings.
\end{theorem}

Probably the reader already knows the reason for one of 
the exceptional cases. We discuss them in \S 2.1. 
The main work begins in \S 2.2, where we show that,
provided $n \ge 7$, there is only one possible way to assign partitions of
the form $(n-m,m)$ to the rows of $X$. (We say that
such partitions are \emph{two-row partitions}.) 
Then in \S 2.3 we show that given such a partial row labelling, 
there is only one way to assign all the column labels.
Of course, once we have fixed the column labels, the remaining row labels
are  uniquely determined. It follows that there is a unique
way to label~$X$. In \S 2.4 we give an efficient way to complete
the row labelling.

It is natural to ask the analogous question for alternating
groups, and for the Brauer character tables of symmetric and
alternating groups. Strikingly, these questions may 
also be answered using the strategy we have just outlined, giving
the results
stated below in Theorem 3.1, Theorem 4.1 and Theorem 5.2
respectively. (We
only deal with alternating groups in odd characteristic.)

Informally, we may interpret our results as saying that the rich structure
of representations of  symmetric and alternating
groups gives their character tables a high degree of rigidity.
In~\S 6 we make this idea more precise by putting 
our results in the  general
context of representations of finite groups.
We also pose two open problems
suggested by our work.

\section{Proof of theorem \ref{Thm:labelling}}

\subsection{}
We begin with $S_6$, the only symmetric group to have an outer
automorphism (see for example \cite[Theorem~7.9]{Rotman}). 
Any outer automorphism of~$S_6$ permutes the conjugacy classes of~$S_6$ by
\[ (6) \leftrightarrow (3,2), \ (3,3) \leftrightarrow (3), \ (2,2,2)
\leftrightarrow (2) \]
and permutes its ordinary characters by
\[ (5,1) \leftrightarrow (2,2,2), \ (2,1,1,1,1) \leftrightarrow (3,3),
\ (4,1,1) \leftrightarrow (3,1,1,1). \]
Thus if $\sigma$ and $\tau$ denote the corresponding permutations
on the set of partitions of $6$ then 
$\chi^\lambda(\nu) = \chi^{\lambda^\sigma}(\nu^\tau)$ for all
partitions $\lambda, \nu$. This gives two different
ways to label the character table of $S_6$. Inspection of the table
shows there are no more.

We now turn to $S_4$, which has the character table shown below.

\medskip
\begin{center}
\begin{tabular}{r|ccccc}
\rule[-2mm]{0mm}{6mm} & $(1^4)$ & $(2,1,1)$ & $(4)$ & $(2,2)$ & $(3,1)$  \\
\hline \rule{0mm}{4mm} $(4)$ & $1$  & $1$ & $1$ & $1$ & $1$ \\
	$(3,1)$ & $3$ & $-1$  & $1$ & $-1$ & $0$ \\
	$(2,1,1)$ &$3$ & $1$ & $-1$ & $-1$ & $0$  \\
	$(2,2)$ & $2$ & $0$ & $0$ & $2$  & $-1$  \\

	$(1^5) $ &$1$ & $-1$ &$1$ &$1$ &$-1$
\end{tabular}
\end{center}
\medskip
Notice that if we swap the columns labelled by $(2,1,1)$ and $(4)$, 
and  the rows labelled by $(3,1)$ and $(2,1,1)$, 
we end up with the same matrix. Again it is easy to see that 
this is the only alternating way to label the character table. 
Unlike the case of $S_6$, this alternative labelling 
is not induced by any automorphism of the group --- its existence 
appears to be entirely coincidental. (See \S 6 for some
remarks related to this phenomenon.)

Theorem~1.1 may readily be verified by inspection of the character tables
if $n \le 3$ or $n=5$, so from now on we shall assume $n \ge 7$.

\subsection{}
Let $n \ge 7$ and let $X$ be an unlabelled character table of $S_n$.
We shall show that there is a unique way to assign two-row partitions
(that is, partitions of the form $(n-m,m)$) to the rows of $X$.

By the orthogonality relations for ordinary characters,
there is a unique row of $X$ containing only positive entries.
Similarly there is a unique column of $X$ containing
only positive entries. We may therefore uniquely identify
the column of the character table corresponding to the
identity element and the row corresponding to the trivial
character, thus fixing the row label $(n)$ and the column label
$(1^n)$. 

For all $n \not= 6$ the symmetric group $S_n$ has exactly
two characters of degree $n-1$, namely $\chi^{(n-1,1)}$
and $\chi^{(2,1^{n-2})}$ (see \cite{JK} Theorem~2.4.10). These
characters are defined by
\[ \chi^{(n-1,1)}(g) = | \Fix g | - 1, \]
where $\Fix g$ is the set of elements fixed by $g$ in
its natural action on $\{1,\ldots,n\}$, and
\[ \chi^{(2,1^{n-2})}(g) = \sgn(g) \chi^{(n-1,1)}(g).\]
For $n > 4$, only $\chi^{(n-1,1)}$
takes $n-3$ as a value. Hence
we may fix the row label $(n-1,1)$.

To proceed further, we observe that
given two rows in the unlabelled character table $X$ we may
multiply the corresponding entries and so obtain a new character.
Then, still working only with the unlabelled table, we may
take the inner product of this character with each
row, and so determine its irreducible constituents. 
The following lemma gives the results needed to exploit these observations. 
(The main idea in the proof of this lemma comes from \cite{Bessenrodt}.)

\begin{lemma}\label{lemma:resind} If $n \ge 4$ then 
\[ \chi^{(n-1,1)}\chi^{(n-1,1)} = \chi^{(n)} + \chi^{(n-1,1)}
+ \chi^{(n-2,2)} + \chi^{(n-2,1,1)}. \]
If $r \ge 2$ and $n > 2r+1$ then
\[
\begin{split} 
\chi^{(n-1,1)} \chi^{(n-r,r)} = \chi^{(n-r-1,r+1)} + \chi^{(n-r,r)} \qquad\qquad
\\ \qquad\qquad+ \chi^{(n-r+1,r-1)} + \chi^{(n-r-1,r,1)} + \chi^{(n-r,r-1,1)}.
\end{split}
\]
\end{lemma}

\begin{proof}
Recall that $\chi^{(n-1,1)} + \chi^{(n)}$ is the permutation
character of $S_n$ acting on $\{1, \ldots, n\}$. As this character
is induced from the subgroup $S_{n-1}$ of $S_n$, we have
\begin{equation*}
\label{resind}
\chi^{(n-1,1)} \theta = \left( 1_{S_{n-1}}\Ind^{S_n} \right)\theta  - \theta
= \theta\Res_{S_{n-1}}\Ind^{S_n} - \theta 
\end{equation*}
for any character $\theta$ of $S_n$. The result now follows from
the branching rule for ordinary representations
of $S_n$ (see~\cite[Ch.~9]{James}).
\end{proof}

This lemma implies that by decomposing the product 
$\chi^{(n-1,1)}\chi^{(n-1,1)}$
we may find the two rows of~$X$ which should be labelled
by $\chi^{(n-2,2)}$ and $\chi^{(n-2,1,1)}$. Since
these irreducible
characters have different degrees (see Lemma~\ref{lemma:dimensions} below)
we may fix the row labels $(n-2,2)$ and $(n-2,1,1)$.

The remaining two-row characters are found in a similar way.
Suppose inductively that there is a unique
way to assign the labels $(n-s,s)$ and $(n-s,s-1,1)$ for $s \le r$.
If $2r = n$ or $2r+ 1 = n$ 
then we are finished, so we assume that $n > 2r+1$.
By decomposing the product $\chi^{(n-1,1)}\chi^{(n-r,r)}$
we may find the two rows which should be labelled
by $\chi^{(n-r-1,r+1)}$ and $\chi^{(n-r-1,r,1)}$. Again, by
Lemma~\ref{lemma:dimensions},
these characters have different degrees, so 
we may fix the row labels
$(n-r-1,r+1)$ and $(n-r-1,r,1)$. This completes the inductive step.

\begin{lemma}\label{lemma:dimensions}
If $r \ge 2$ and $n \ge 2r$ and  then
\[ \chi^{(n-r,r)}(1) < \chi^{(n-r,r-1,1)}(1). \]
\end{lemma}

\begin{proof}
By the hook-formula for the irreducible
character degrees of $S_n$ (see \cite[Theorem~2.3.21]{JK}),
\[
\chi^{(n-r,r)}(1) = \binom{n}{r} \frac{n-2r+1}{n-r+1} \]
which is  less than
\[ \chi^{(n-r,r-1,1)}(1) = (r-1) \binom{n}{r}  \frac{n-2r+2}{n-r+2}. \]

\vspace{-22pt}
\end{proof}
\phantom{D}
\vspace{-12pt}

\subsection{}
Fix a column of $X$. Suppose that in one possible
labelling of the columns of~$X$, the corresponding conjugacy class 
has~$g\in S_n$
as a representative.
We shall
show that the cycle type of $g$ can be reconstructed
from the character values $\chi^{(n-r,r)}(g)$, which are known
from the previous section. 
This will show that the column labels of $X$ are uniquely determined.

Since $\chi^{(n-1,1)}(g) = 
|\Fix g| - 1$, we can find the number of fixed points
of~$g$. 
Suppose inductively that
we know that $g$ has $a_1$ fixed points, $a_2$ $2$-cycles,
\ldots, and $a_{r-1}$ $(r-1)$-cycles, where
$1 < r \le n/2$. Let
\[ \pi^{(n-r,r)}=  \sum_{s=0}^{r}
\chi^{(n-s,s)}. \]
This is the permutation character of $S_n$ acting on $r$-subsets of
$\{1, \ldots, n\}$. Using only the known character values 
we may calculate $\pi^{(n-r,r)}(g)$, and so
find the number of $r$-subsets fixed by~$g$. This
equals the number of $r$-cycles in~$g$,
plus some further quantity
that can be computed given $a_1, \ldots, a_{r-1}$. Therefore~$a_r$ may be determined. Since the action of $S_n$ 
on $r$-subsets is isomorphic to its action on $(n-r)$-subsets,
this is sufficient to determine the number of cycles in~$g$
of any given length.

\subsection{}
Now we have found all the column labels of $X$, the remaining
row labels are uniquely determined. 
Here
we briefly give a practical way to determine these labels.

It follows from a result of Kramer \cite{Kramer}
that two irreducible
characters of $S_n$
which agree on every cycle of length $r$ for $1 \le r \le n$
are equal. (See \cite{WildonCharValues}
for a short proof of this, and some related results.)
Suppose then that we have already computed the character values
$\chi^{\lambda}(n-r,1^r)$ where $0 \le r < n$ and
$\lambda$ is any partition of $n$. (One
convenient way to do this is to use
the Murnagham--Nakayama rule and the hook-formula.)
Then to find the partition labelling row $i$ of the 
character table $X$, we need only compare the
values $X_{i\nu}$ for
$\nu = (n-r,1^r)$ and $0 \le r < n$ with
those in our pre-computed table.

\section{Alternating groups} 
The outer automorphisms of the alternating group $A_n$
induced by the conjugacy action of $S_n$
leads to some inevitable ambiguity
in the row and column labels of its character table.
Recall that
any such automorphism acts as an involution, swapping pairs of  split
characters and conjugacy classes.
More precisely, if $\lambda$ is a partition of $n$
then $\chi^\lambda\res_{A_n}$ is reducible if and only
if $\lambda$ is self-conjugate, in which case it splits
as a sum of two irreducible characters of $A_n$. Similarly the
conjugacy class labelled by the partition $\nu$ of $n$ 
splits in $A_n$ if and only
if $\nu$ has odd distinct parts. We distinguish
split characters and classes by arbitrarily allocating~$+$ and~$-$ signs.
Whichever allocation we choose, it will be reversed under 
the outer action of $S_n$.

For
example, the character table of $A_5$ is

\medskip
\begin{center}
\begin{tabular}{r|rrrrr}
\rule[-3mm]{0mm}{8mm} & $(1^5)$ & $(2,2,1)$ & $(3,1,1)$ & $(5)^+$ & $(5)^-$  \\
\hline \rule{0mm}{4mm} $(5)$ & $1$ & $1$& $1$& $1$ & $1$\\
	$(4,1)$ & $4$ & $0$& $1$ & $-1$ & $-1$ \\ 
	$(3,2)$ & $5$ & $1$ & $-1$ & $0$ & $0$ \\ 
$(3,1,1)^+$ & $3$ & $-1$ & $0$ & $\alpha$ & $\beta$ \\ 
$(3,1,1)^-$ & $3$ & $-1$ & $0$  & $\beta$ & $\alpha$ \\	
\end{tabular}
\end{center}
\medskip

\noindent where $\alpha = (1+\sqrt{5})/2$ and $\beta=(1-\sqrt{5}/2)$.
The alternative labelling induced by the action of $S_5$ is 
given by swapping
the characters and conjugacy classes
labelled~$+$ and $-$.

The next theorem states that this is usually the only alternative
labelling.

\begin{theorem}
Let $X$ be a character table of the alternating group $A_n$. 
Provided $n \not= 6$
there is a unique way to assign non-self-conjugate partitions
to the rows of $X$ and partitions not having odd distinct parts
to the columns of 
$X$ so that $X_{\lambda\nu} = \chi^{\lambda}(\nu)$ for
all such $\lambda$ and $\nu$.
The labels of the split characters and conjugacy classes
are uniquely determined up to signs.
\end{theorem}

\begin{proof} For $n \le 5$ the theorem can be readily
verified by inspecting the tables. 
By \cite[Theorem 2.5.15]{JK},
provided
$n\ge 7$, the only character of $A_n$ of degree $n-1$ is
the one labelled by $(n-1,1)$. So we can fix this label.
Since
no self-conjugate partitions appear in the calculations of \S 2.2,
the remaining two row characters may then be identified as before. 
The cycle types of the columns may also be identified as before.
When one comes to a split class the labels $+$ and $-$ may be assigned
either way round. Once all the column labels
are fixed, there is then a unique way to label
the rows.
\end{proof}

Inspection of the character table of $A_6$ shows that it has $4$ alternative
labellings. Using the definitions given in \S 6 below, we have 
$\clAut(A_6) \cong \cAut(A_6) \cong \Out(A_6) \cong \left<(12),(34)\right>$.

\section{Brauer character tables of symmetric groups}

Let $F$ be a field of prime characteristic $p$. We say
that a partition
is $p$-regular if it has at most $p-1$ parts
of any given size.
Recall from \cite{James} that the irreducible
$FS_n$ representations are parametrised
by the $p$-regular partitions of $n$. Let $D^\mu$
be the $p$-modular irreducible corresponding
to the $p$-regular partition $\mu$ and let~$\phi^\mu$ be the 
Brauer character of $D^\lambda$.
(See \cite{Navarro} for an introduction to Brauer characters.)
If $\nu$ is
a partition of $n$ with no part divisible by $p$ then we write~$\phi^\mu(\nu)$ for the value
of~$\phi^\mu$ on the conjugacy class of $p'$-elements labelled
by~$\nu$. 

In this section we shall prove the following theorem.

\begin{theorem}\label{thm:modp}
Let $X$ be a Brauer character table of the symmetric group~$S_n$
in characteristic $p$. Unless $n = 6$, or $n=4$ and $p > 2$,
there is a unique way to assign $p$-regular partitions
to the rows of $X$ and partitions with no part divisible by $p$
to the columns of $X$ so that $X_{\mu\nu} = \phi^{\mu}(\nu)$ for
all such $\mu$ and~$\nu$. In the exceptional cases there
are exactly $2$ different labellings.
\end{theorem}

To prove this theorem
we shall need some further results on the modular representations of symmetric
groups. Recall that the \emph{decomposition matrix}~$D_n(p)$ of $S_n$
in characteristic $p$ is defined by
\[ \chi^\lambda = \sum_{\mu} D_n(p)_{\lambda\mu} \phi^\mu \]
when $\lambda$ is a partition of $n$ and the sum is over
all $p$-regular partitions $\mu$.
(Here and elsewhere when we write a relation between ordinary
and Brauer characters it is intended to hold 
for $p'$-elements only.) 
Let~$\unrhd$ denote the dominance order
on partitions of $n$ (see \cite[Definition~3.2]{James}).
The following lemma is Corollary 12.3 in \cite{James}.

\begin{lemma}\label{lemma:wedge}
Let $\lambda$ be a partition of~$n$ and let
$\mu$ be a $p$-regular partition of~$n$. If
$D_p(n)_{\lambda\mu} \not=0$ then $\mu \unrhd \lambda$.
Moreover $D_p(n)_{\lambda\lambda} = 1$.\hfill$\Box$
\end{lemma}

We also need a simple branching rule for modular
representations. 

\begin{lemma}\label{lemma:res}Let $\mu = (\mu_1,\ldots,\mu_k)$ be a partition of $n$ such
that $\mu_1 > \mu_2$. If $\bar{\mu} = (\hbox{$\mu_1 -1$},
\mu_2,\ldots, \mu_k )$ is $p$-regular then 
$\phi^\mu\res_{S_{n-1}} = \phi^{\bar{\mu}} + \psi$
where $\psi$ is a sum of Brauer characters $\phi^\lambda$
labelled by partitions $\lambda$ such that $\lambda \rhd \bar{\mu}$.
\end{lemma}

\begin{proof}In future we shall write $\psi \rhd \mu$ if $\psi$
is a weighted sum of Brauer (or ordinary) characters labelled
by partitions $\lambda$ such that $\lambda \rhd \mu$.
By Lemma~\ref{lemma:wedge} we may write
$\phi^\mu = \chi^\mu - \theta$ where $\theta \rhd \mu$.
It then follows from the ordinary branching rule that
\[ \phi^\mu \Res_{S_{n-1}} = \chi^{\bar{\mu}} + \psi \]
where $\psi \rhd \bar{\mu}$. Now apply
Lemma~\ref{lemma:wedge} one more time.
\end{proof}

We can now begin the proof of Theorem~\ref{thm:modp}.
We follow  as closely as possible the method of proof used in \S 2; thus
\S 4.1 below is the analogue of~\S 2.1, and so on. 

\subsection{} If $n \le 6$ the theorem may readily be verified 
by inspecting the tables. (Brauer character tables for $p=2$ and
$p=3$ and $n\le 10$ appear in Appendix~I.F
of~\cite{JK}; for $p=5$ and $n=5,6$ the required tables may  easily be
calculated by hand, as all blocks have weight $0$ or $1$.)
The only difference in behaviour from the ordinary case  occurs when $n=4$: the
Brauer character table of $S_4$ in characteristic $2$ is
shown below.

\medskip
\begin{center}
\begin{tabular}{r|cc}
\rule[-2mm]{0mm}{6mm} & $(1^4)$ & $(3,1)$   \\
\hline \rule{0mm}{4mm} $(4)$ & $1$  & $1$  \\
	$(3,1)$ & $2$ & $-1$  
\end{tabular}
\end{center}
\medskip

\noindent Clearly there is no longer any ambiguity about the labels.


\subsection{}Let $n \ge 7$ and let $X$ be an unlabelled Brauer
character table of $S_n$ in characteristic $p$.
As before, it is easy to see that there is a unique way
to assign the row label $(n)$ and the column label $(1^n)$.
It was first proved by Wagner (see \cite{WagnerEven,WagnerOdd})
that
if~$S$ is a simple $FS_n$-module with $\dim S \le n-1$
then, provided $n \ge 7$, either $S$ is $1$-dimensional, or~$S$ is isomorphic
to one of $D^{(n-1,1)}$
or $D^{(n-1,1)} \otimes \sgn$. (For an alternative shorter proof
see James~\hbox{\cite[Theorem 6]{JamesMinimal}}.)
The Brauer character of $D^{(n-1,1)}$ is  
\[ \phi^{(n-1,1)}(g) = |\Fix g| - \begin{cases} 1 & \text{if $p \notdiv n$} 
\\ 2 & \text{if $p \mid n$}. \end{cases}\]
Of $\phi^{(n-1,1)}$ and $\phi^{(n-1,1)} \sgn$,
only $\phi^{(n-1,1)}$ takes $n-4$ as a value,
hence we may identify the row of $X$ labelled by $(n-1,1)$.

We are now in a position to  identify the rows
of $X$ labelled by all two-row partitions. As before, we
do this inductively by taking tensor products.
However, as there
is no simple formula for
the degrees of the characters~$\phi^{(n-r,r)}$, we have
to be slightly more subtle in our approach. 
The following lemma is the
analogue of Lemma~\ref{lemma:resind}.

\begin{lemma}\label{lemma:modprod} 
If $r \ge 1$ and $n > 2r+1$ then 
\[ \phi^{(n-1,1)}\phi^{(n-r,r)} = b\phi^{(n-r-1,r+1)} + \phi^{(n-r-1,r,1)} 
+ \psi \]
where $b \ge 1$ and $\psi$ is a weighted sum of irreducible Brauer
characters labelled by partitions $\mu$ such that $\mu \unrhd (n-r,r-1,1)$.
\end{lemma}

\begin{proof}
We adapt the argument used in Lemma~\ref{lemma:resind}.
Since $\phi^{(n-1,1)} = 1_{S_n-1}\ind^{S_n} - c\chi^{(n)}$ 
for some $c \in \{1,2\}$,
it
is sufficient to show that
\[ \phi^{(n-r,r)}\Res_{S_{n-1}}\Ind^{S_{n}} = \phi^{(n-r-1,r+1)}
+ \phi^{(n-r-1,r,1)} + \psi
\]
where $\psi \unrhd (n-r,r-1,1)$. We have
\[ \phi^{(n-r,r)} = \chi^{(n-r,r)} - \theta \]
where $\theta \unrhd (n-r+1,r-1)$. Hence by the
the ordinary branching rule,
\[ \phi^{(n-r,r)}\Res_{S_{n-1}}\Ind^{S_{n}} =
\chi^{(n-r-1,r+1)} + \chi^{(n-r-1,r,1)} + \psi\]
where $\psi \unrhd (n-r,r-1,1)$. Now apply
Lemma~\ref{lemma:wedge}.
\end{proof}

For simplicity we state the following proposition for~$p > 2$, 
and explain later the
small modifications needed if~$p=2$.

\begin{proposition}\label{prop}
Let $n \ge 7$ and 
let $C$ be a labelled Brauer character table of~$S_n$ in
characteristic~$p > 2$. Suppose that the first
column is labelled by~$(1^n)$, and that the rows are arranged 
so that the first row is labelled by~$(n)$, the second by $(n-1,1)$
and the next $2(r-1)$ with the labels
\begin{align*}  &(n-2,2), (n-3,3), \ldots, (n-r,r) \\
&(n-2,1^2), (n-3,2,1), \ldots, (n-r,r-1,1).\end{align*}
in \emph{any order}. 
If we are given the first $2r$ rows of $C$ with the row
and column labels removed, then the row labels may
be uniquely reconstructed.
\end{proposition}

\begin{proof}
We work by induction on $n$. If $n = 7$ and $p=3$ or $p=5$
then the Brauer characters
of~$S_n$   that can appear in~$X$ 
have distinct degrees,
so the result is immediate. If $p=7$ then both $\phi^{(5,2)}$ and 
$\phi^{(4,3)}$ have degree $14$, but only the former
takes $6$ as a value. (As both characters are in blocks of weight zero
this can be seen directly from the ordinary character
table.)

Suppose now that $n \ge 8$.
By hypothesis, we may immediately attach the row labels $(n)$
and $(n-1,1)$ to $C$.
We now attempt to reach a situation in which the inductive
hypothesis for $n-1$ can be applied. Notice first that 
the values of $\phi^{(n-1,1)}$ determine which
columns in the table come from conjugacy 
classes with at least one fixed point,
and so are relevant when we restrict a character to~$S_{n-1}$.
The restriction of~$\phi^{(n)}$ to~$S_{n-1}$ is of course~$\phi^{(n-1)}$.
We may obtain~$\phi^{(n-2,1)}$ by removing
any copies of~$\phi^{(n-1)}$ from $\phi^{(n-1,1)}\res_{S_{n-1}}$.

By Lemma~\ref{lemma:modprod} above, when we express
$\phi^{(n-1,1)}\phi^{(n-1,1)}$ as a sum of rows of~$X$, two
new characters appear: $\phi^{(n-2,2)}$ and $\phi^{(n-1,1,1)}$.
When we restrict these new characters to $S_{n-1}$ we get,
in addition to any copies of $\phi^{(n-2,1)}$ and $\phi^{(n-2)}$
that may be present, two new Brauer characters of 
$S_{n-1}$: namely~$\phi^{(n-3,1,1)}$ and~$\phi^{(n-3,2)}$. 
We may now apply the inductive hypothesis (with~$r=2$) to
determine which label should go with which. To get back
to~$S_n$ we use~Lemma~\ref{lemma:res}. Together
with Lemma~\ref{lemma:wedge}, it implies that
$\phi^{(n-2,2)}\res_{S_{n-1}}$ does not contain~$\phi^{(n-3,1,1)}$,
whereas $\phi^{(n-2,1,1)}\res_{S_{n-1}}$ does. We use this to
fix the labels~$(n-2,2)$ and~$(n-2,1,1)$.

The remaining two-row labels are fixed by repeating
this argument, in a way closely
analogous to the proof of Lemma~2.2. We therefore leave the remaining details of the proof to the reader.
\end{proof}

By a further induction on $r$ we may use this proposition
to identify the rows of $X$ labelled by two-row partitions.

If $p=2$ then the statement of Proposition~\ref{prop} must
be slightly modified. If $n$ is odd we must delete $(n-2,1,1)$,
and if $n = 2m$ is even we must delete $(n-2,1,1)$ and $(m,m)$, as 
these partitions are no longer $p$-regular.
The main change
in the proof is that now $\phi^{(n-1,1)}\phi^{(n-1,1)}$ only
contains one new Brauer character, $\phi^{(n-2,2)}$; this makes the
first step is slightly simpler. After that, no alterations are needed,
unless $n=2m$ is even, in which case the last two
row Brauer character we must find is $\phi^{(m+1,m-1)}$.
Again if anything this makes the process slightly simpler.

\subsection{}
We now determine the column labels of $X$. If we order
the rows and columns of the decomposition matrix $D_n(p)$ by
the dominance order, $(n)$ appears first of all, followed
by the two-row partitions of $n$. By Lemma~\ref{lemma:wedge},
$D_n(p)_{\lambda\mu} = 0$
if~$\lambda$ has at most two rows and $\mu$ does not. Hence
the matrix~$D_n(p)X$ has at its top 
the values of the characters~$\chi^\lambda$
for partitions~$\lambda$ with at most two rows.

We can now use the same argument as in the ordinary case
to determine the column labels of $X$. Once we have
a complete set of column labels the row labels are, of course,
fixed. This completes the proof of Theorem 4.1.

\section{Brauer character tables of alternating groups}
Provided $p$ is odd, the modifications 
to the work of \S 4 needed
to deal with the Brauer character tables of alternating groups
in characteristic $p$ are analogous to those needed
to the work of \S 2 (and made in \S 3) to
deal with the ordinary character
tables of alternating groups. 

Let $X$ be a Brauer character table of the alternating group $A_n$
in odd-characteristic.
As in \S 3, we label the split
characters and conjugacy classes by~$+$ and~$-$ signs.
As usual, to get started we need to identify the character labelled
by \hbox{$(n-1,1)$}. For this we use Theorem 1.1 in \cite{WagnerOdd},
which states
that if $n\ge 7$ and
$\phi$ is an odd-characteristic Brauer character of~$A_n$
such that $\phi(1) \le n$ then $\phi = \phi^{(n-1,1)}\res_{A_n}$.
(For an alternative shorter proof for $n \ge 10$ see 
\cite[Theorem~7(ii)]{JamesMinimal}; the result can easily be checked directly
in the remaining cases.)

We also need to know that none of the characters considered in 
\S 4.2 split on restriction to $A_n$. For this to hold
we must take $n \ge 8$.

\begin{lemma}If  $\lambda$ is a $p$-regular partition of $n$
let $m(\lambda)$
be the $p$-regular partition of $n$ defined by
\[ D^\lambda \otimes \sgn = D^{m(\lambda)}.\]
If either $\lambda$ has at most two rows, or $\lambda$ is of the form
$(n-r,r-1,1)$ for some $n \ge 8$ and $r \le n/2$
then  $\lambda \not= m(\lambda)$. Hence
$\phi^\lambda$ does not split on restriction to $A_n$.
\end{lemma}

\begin{proof}It
follows easily from Ford's description of the
Mullineux map in~\cite{FordIrreps}
that $\lambda \not= m(\lambda)$ if either
of the conditions on~$\lambda$ hold.
Ford's paper
also gives the Clifford theory needed to prove the second assertion.
\end{proof}

It is not possible to take $n \ge 7$ (as was the case in \S 4.2) 
because if $p=3$
then $m((4,2,1))= (4,2,1)$.
The base case in the analogue of Proposition~\ref{prop}
is therefore $n=8$. Calculation shows that the Brauer characters
of $A_8$ that can appear in the table $C$ have distinct
degrees when $p=3$ and when $p=7$. If $p=5$ then $\phi^{(6,1,1)}
= \phi^{(5,3)} = 21$, but only the former character takes $6$ as a value,
so again there is no ambiguity. Thus we may identify
the rows of $X$ labelled by two-row partitions. The column
labels may now be determined in essentially the same way
as \S 4.3.

Direct examination of the cases for $n \le 7$ gives
the following theorem.

\begin{theorem}
Let $X$ be a Brauer character table of the alternating group~$A_n$
in odd characteristic $p$. 
Provided $n \not= 6$
there is a unique way to assign $p$-regular partitions
$\lambda$ such that $m(\lambda) \not= \lambda$ 
to the rows of $X$ and partitions not all of whose
parts are odd, and with no part divisible by $p$, to the columns of 
$X$ so that $X_{\lambda\mu} = \phi^{\lambda}(\mu)$ for
all such $\lambda$ and $\mu$.
The labels of the split characters and conjugacy classes
are uniquely determined up to signs.\hfill$\Box$
\end{theorem}

When $n=6$ and $p=3$ there are two different labellings, interchanged
by the conjugacy action of $S_6$. When $n=6$ and $p=5$
there are again two different labellings, but this time $S_6$
acts trivially, and they
are interchanged only by the outer automorphism of $S_6$.

\section{A more general setting}

Given a arbitrary finite group $G$, there is usually no canonical
way to label the rows and columns of its character table. However, as 
the
reader familiar with Brauer's permutation lemma will
already have realised, this need not stop
us from considering analogous versions of our results.

Let $k$ be the number of conjugacy classes of $G$, and let
$X$ be a character table of $G$. 
We say that a pair $(\sigma,\tau) \in S_k \times S_k$ 
is an \emph{automorphism} of~$X$ if 
$X_{i\sigma j\tau} = X_{ij}$ whenever $1 \le i,j \le k$. 
Let $\Aut(X)$ be the group of
all automorphisms of $X$. 
It is clear that for each $\sigma \in S_k$
there is at most one $\tau \in S_k$ such that $(\sigma,\tau) \in \Aut(X)$.
We may therefore define a group~$\cAut(G)$ by
\[ \cAut(G) = \{ \sigma \in S_k : \text{$(\sigma,\tau) \in \Aut(X)$
for some $\tau \in S_k$} \}. \]
This group is well-defined up to conjugacy in $S_k$. 
For example, our Theorem~1.1 states that $\cAut(S_n)$
is trivial unless $n=4$ or $n=6$.

\begin{problem}
Calculate $\cAut(G)$ for  important classes of groups.
\end{problem}

In connection with this problem, it is interesting to explore the 
relationship between 
$\cAut(G)$ and $\Out(G) = \Aut(G)/\Inn(G)$, the group of 
outer automorphisms of $G$. Clearly
there is a group homomorphism
\[ c : \Out(G) \rightarrow \cAut(G) \]
defined for $\gamma \in \Out(G)$ by letting $c(\gamma)$ be the permutation 
induced by $\gamma$ on the ordinary characters of $G$. In some
cases $c$ is an isomorphism --- for example, this is the case
if $G$ is abelian, or $G$ is a symmetric group other than $S_4$,
or $G$ is any alternating group.
But, as the example of $S_4$ shows, $c$ need not be surjective.
The dual question, of whether $c$ must be injective,
or equivalently, whether there
is a finite group~$G$ and an outer automorphism $\gamma \in \Aut(G)$ 
such
that~$\chi^\gamma = \chi$ for all irreducible characters~$\chi$,
was considered by Burnside: see Note B in~\cite{Burnside}. 
In \cite{WallOuter}, 
G.~E.~Wall gives
an example in which~$G$ has order $32$ and $\gamma$ has order $4$.

It is worth noting that Burnside's question can be stated 
without even mentioning characters,
since by Brauer's permutation lemma (see \cite[\S 6]{BrauerPL}),~$t$ 
is such an automorphism
if and only if $t$ permutes within themselves all the conjugacy classes
of $G$. (Incidentally, it seems clear from \cite[\S 186]{Burnside} 
that Brauer's permutation lemma was already well known to Burnside.) 

Another obvious question, which is related to Problem 6.1, is:

\begin{problem} Is the map $c : \Out(G) \rightarrow \cAut(G)$
always an isomorphism when $G$ is a finite simple group?
\end{problem}

A final problem, which can be answered more easily, 
arises from the definition of $\cAut(G)$.
If we look instead at the admissible \emph{column} permutations of $X$ 
then we obtain the group
\[ \clAut(G) = \{ \tau \in S_k : \text{$(\sigma,\tau) \in \Aut(X)$
for some $\sigma \in S_k$} \}. \]
By Brauer's permutation lemma  
the groups $\cAut(G), \clAut(G) \le S_k$ are isomorphic as
abstract groups, \emph{via} an isomorphism preserving the cycle
types of elements. But this on its own does not guarantee
that they are permutation isomorphic, 
as the two subgroups of $S_6$,
\[ \left<(12)(34), (13)(24)\right>, \quad \left<(12)(34), (12)(56)\right> \]
show. (There are many more examples of this type.) The following
example shows that
$\cAut(G)$ and $\clAut(G)$ need \emph{not} be permutation isomorphic,
and so we made a genuine choice in concentrating on $\cAut$ earlier.

\begin{example}Let $G \cong C_2 \times D_8$, where $D_8$ is the dihedral
group of order~$8$. The character table of $G$ is, with one
ordering of the rows and columns: 
{\small
\[ \begin{matrix}
-1  & 1 & 1 &-1 & 1 &-1 &-1 &-1 & 1 & 1 \\
-1 & -1 & 1 & 1 &-1 & 1 &-1 &-1 & 1 & 1 \\
 1 &-1& -1&  1 & 1 &-1& -1& -1  &1 & 1 \\
 1 & 1& -1& -1 &-1&  1& -1& -1&  1&  1 \\
-1&  1& -1&  1& -1& -1&  1&  1&  1&  1 \\
 1& -1&  1& -1& -1& -1&  1&  1&  1&  1 \\
-1& -1& -1& -1&  1&  1&  1&  1&  1&  1\\
 0&  0&  0&  0&  0&  0 &-2&  2&  2& -2 \\
 0 & 0&  0&  0&  0&  0&  2& -2&  2& -2 \\
 1 & 1&  1 & 1&  1&  1&  1&  1&  1&  1 
\end{matrix} \]}
One finds that
\begin{align*} \cAut(G) &= \left< (1234)(56), (34)(57), (89) \right>, \\
 \clAut(G) &= \left< (1234)(56), (25)(46), (78) \right>
\end{align*}
where the isomorphism $\cAut(G) \cong \clAut(G)$ is indicated
by the order of generators. As the orbits of $\cAut(G)$
have sizes $4,3,2,1$ whereas the orbits of $\clAut(G)$ have
sizes $6,2,1,1$, the two groups are not permutation
isomorphic. (Abstractly, both  are isomorphic to $S_4 \times C_2$.)
\end{example}

Finally we  mention
a theorem of Higman (see 
\cite[Theorem~8.21]{IsaacsChars}) which
states that given a character table of a finite group, 
one can determine prime divisors of the orders of the group
elements corresponding to any given column.
It is well known that
the dihedral and quaternion groups of order $8$ have the same character table, so this is the most one can hope for in general. Theorem~1.1 and
Theorem~3.1 imply
that for symmetric and alternating groups much more is true.

\begin{corollary}
Given an unlabelled character table
of a symmetric group other than $S_4$ 
one may determine the order of the elements corresponding to
any of its columns. The same result holds for any alternating group. 
\hfill$\Box$
\end{corollary}

\section*{Acknowledgement}
I should like to thank Dr Christine Bessenrodt for the helpful
comments she made when this work was beginning.

\def\cprime{$'$} \def\Dbar{\leavevmode\lower.6ex\hbox to 0pt{\hskip-.23ex
  \accent"16\hss}D}

\end{document}